\documentclass[fontsize=12pt,a4paper,headings=normal,
twoside=false,leqno,parskip=half-,abstract=true]{scrartcl}
\usepackage[english]{babel}
\usepackage[utf8]{inputenc}
\usepackage{hyperref}
\hypersetup{
 pdftitle={Lyapunov function},
 pdfauthor={Phillipo Lappicy},
 colorlinks=true,
 linkcolor=blue,
 citecolor=blue,
 filecolor=blue,
 urlcolor=blue}

\usepackage{graphicx}
\usepackage[format=plain,labelfont=bf,font=small]{caption}
\usepackage{xcolor,colortbl}
\definecolor{Gray}{gray}{0.9}
\usepackage[arrow, matrix, curve]{xy}
\usepackage{float}

\usepackage{caption}
\usepackage{tabulary}
\usepackage{array}
\newcolumntype{N}[1]{>{\centering\arraybackslash}m{#1}}

\usepackage{amsmath,amsthm}
\usepackage{amssymb}
\usepackage{bigints}
\usepackage{comment}

\makeatletter
\newcommand{\tpitchfork}{%
  \vbox{
    \baselineskip\z@skip
    \lineskip-.52ex
    \lineskiplimit\maxdimen
    \m@th
    \ialign{##\crcr\hidewidth\smash{$-$}\hidewidth\crcr$\pitchfork$\crcr}
  }%
}
\makeatother
\usepackage{latexsym}

\usepackage[notref,notcite,color,
final
]{showkeys}
\usepackage[pagewise]{lineno}

\usepackage{tikz}

\usepackage[textwidth=2cm,textsize=small,backgroundcolor=none]{todonotes}

  \mathchardef\ordinarycolon\mathcode`\:
  \mathcode`\:=\string"8000
  \begingroup \catcode`\:=\active
    \gdef:{\mathrel{\mathop\ordinarycolon}}
  \endgroup

\theoremstyle{plain}
\newtheorem{thm}{Theorem}[section]

\begin{document}

\title{{\Large{An energy formula for fully nonlinear degenerate parabolic equations in one spatial dimension}}}

\author{
 \\
{~}\\
Phillipo Lappicy*$^\circ$ and Ester Beatriz*\\
\vspace{2cm}}

\date{ }
\maketitle
\thispagestyle{empty}

\vfill

$\ast$\\
ICMC, Universidade de S\~ao Paulo\\
Av. trabalhador são-carlense 400, 13566-590, São Carlos, Brazil\\
$\circ$\\
Universidad Complutense de Madrid\\
Pl. de las Ciencias, 3, 28040, Madrid, Spain\\


\newpage
\pagestyle{plain}
\pagenumbering{arabic}
\setcounter{page}{1}

\begin{abstract}
\noindent Energy (or Lyapunov) functions are used to prove stability of equilibria, or to indicate a gradient-like structure of a dynamical system. Matano constructed a Lyapunov function for quasilinear non-degenerate parabolic equations. 
We modify Matano's method to construct an energy formula for fully nonlinear degenerate parabolic equations. 
We provide several examples of formulae, and in particular, a new energy candidate for the porous medium equation.



\end{abstract}

\section{Main results}

\numberwithin{equation}{section}
\numberwithin{figure}{section}
\numberwithin{table}{section}

We consider the  scalar fully nonlinear partial differential equation
\begin{equation}\label{FULLYEQ}
    f(x,u,u_x,u_{xx},u_t)=0,
\end{equation}
for $x\in (0,1)$ and $t>0$ \textcolor{black}{with appropriate initial data $u_0(x)$}. Here indices abbreviate partial derivatives. We assume that $f\in C^2$ satisfies the following degenerate parabolic conditions 
\begin{equation}\label{par}     
  f_q \cdot f_r \leq 0,  \qquad \text{ and } \qquad f_r \neq 0,
\end{equation}
for every argument $(x,u,p,q,r):=(x,u,u_x,u_{xx},u_t)\in [0,1]\times \mathbb{R}^4$. 
{\color{black}
Conditions \eqref{par} imply that only the diffusion coefficient $f_q$ may vanish, since $f_r \neq 0$ excludes time-evolution type degeneracies\footnote{Note time-evolution degeneracies can be transformed into singular diffusion, 
see \cite[Problem 3.6]{Vazquez06}.} such as in \emph{Trudinger's equation}, $(u^{\alpha})_t=u_{xx}$, for $\alpha>0$; see \cite{Trudinger}. 
Without loss of generality, we consider $f_r < 0$ and thus $f_q \geq 0$. Indeed, if $f_r > 0$, then $f_{\tilde{r}}<0$ for $\tilde{r} := -r$. 
\textcolor{black}{Moreover, in order to guarantee that the diffusion $f_q$ degenerates in a meagre set, we also assume that the following set is of (Lebesgue) measure zero in $[0,1]\times \mathbb{R}^2$,}
\begin{equation}\label{pardeg}     
    \left\{ (x,u,p) \in [0,1]\times \mathbb{R}^2 \text{ $|$ } f_q(x,u,p,0,0)=0 \right\}.
\end{equation}
In particular, the condition \eqref{pardeg} prevents 
that $f_q(x,u,p,0,0)= 0$ for all $(x,u,p)\in [0,1]\times \mathbb{R}^2$ and degeneracies of the same order as $u_{xx}$, such as the \emph{dual porous medium equation}, $u_t=|u_{xx}|^{m-1}u_{xx}$, for $m>1$, see \cite{Bernis,Vazquez06}.
}

We consider \eqref{FULLYEQ} with two types of \emph{separated boundary conditions} at $x=\iota\in \{ 0,1\}$. For each boundary point $x=\iota$, separately, we either assume homogeneous Dirichlet boundary conditions or nonlinear boundary conditions of Robin type, respectively 
\begin{subequations}\label{bc}
	\begin{align}
        u&=0,\label{Dir}\\
        u_x&=b^\iota(u).\label{Rob} 
	\end{align}
\end{subequations}
We assume $b^\iota\in C^1$. Neumann boundary conditions occur if $b^\iota(u)=0$. 
See \cite{Amann88,Lunardi95} for abstract settings involving nonlinear boundary conditions of type \eqref{Rob}.

Equations \eqref{FULLYEQ}-\eqref{bc} include classical examples, such as evolution involving $p$-laplacian diffusion, the porous medium equation or certain mean curvature flow. 
These classical equations with further nonlinear gradient-dependent forcing did not have any apparent variational structure, which we are now able to display. 
It is the scope of this paper to provide a unifying variational formulation to several degenerate fully nonlinear parabolic equations in one spatial dimension. 

Below we construct a \emph{Lyapunov function}
\begin{equation}\label{IntroLyap}     
    E:=\int_{0}^{1} L(x,u,u_x)\,dx  \qquad \text{ such that } \qquad \frac{dE}{dt} < 0
\end{equation}
along non-equilibrium solutions $u=u(t,x)$ of \eqref{FULLYEQ}. 
Therefore, the time-dependent energy $t\mapsto E(u(t,.))$ decreases strictly, except at equilibria, i.e. $u_t\equiv 0$.

Before we present the main result, 
we rewrite the fully nonlinear equation \eqref{FULLYEQ} suitably, following the spirit of \cite{LappicyFiedler19}. Then we modify of Matano's original idea in \cite{Matano88} 
in order to incorporate degeneracies of the PDE \eqref{FULLYEQ} for a Lyapunov function $E$ as in \eqref{IntroLyap}.

Indeed, we split the equation \eqref{FULLYEQ} in order to emphasize the \emph{degenerate diffusion},
\begin{equation}\label{FULLYDIFF}
    F(x,u,u_x,u_{xx},u_t)= f_q(x,u,p,0,0)u_{xx},
\end{equation}
where $F(x,u,p,q,r) := - f(x,u,p,q,r) + f_q(x,u,p,0,0)u_{xx}$. 
The degeneracy conditions in equation \eqref{par} become
\begin{equation}\label{par2}     
    F_r(x,u,p,q,r) > 0, \quad \text{ and }\quad f_q(x,u,p,0,0) \geq F_q(x,u,p,q,r). 
\end{equation}
for every argument $(x,u,p,q,r):=(x,u,u_x,u_{xx},u_t)$. 
\textcolor{black}{Condition \eqref{pardeg} implies that $F(x,u,p,q,r)=-f(x,u,p,q,r)$ for a set of measure zero.} 

Next, we split $F$ into two parts: one which is independent of $u_{xx}$ and $u_t$, whereas another that depends on them. 
First, we distinguish a term $F^0$ related to \emph{reaction}, when $u_{xx}=u_t=0$. 
Second we describe the \emph{time evolution} term $F^1$, the only term that depends on $u_t$. 
Specifically, we define
\begin{align}\label{F0}
\begin{split}
    F^0(x,u,p)&:=F(x,u,p,0,0),\\
    F^1(x,u,p,q,r)&:= F(x,u,p,q,r) - F^0(x,u,p),
\end{split}
\end{align}
where $F^0\in C^2$ and $F^1 \in C^1$, since $f\in C^2$.

The parabolic equation \eqref{FULLYDIFF} can be rewritten as 
\begin{equation}\label{diffEQ}
F^1(x,u,p,q,r) = f_q(x,u,p,0,0)u_{xx} - F^0(x,u,p).
\end{equation}
The degeneracy conditions \eqref{par} incarnated in \eqref{par2} imply that
\begin{equation}\label{par3}     
    F^1_r > 0 \quad \text{and} \quad f_q(x,u,p,0,0) \geq F_q^1(x,u,p,q,r). 
\end{equation}
for every $(x,u,p,q,r):=(x,u,u_x,u_{xx},u_t)$. 
%

The main modification of Matano's method is a different Ansatz for the function $L$ in \eqref{IntroLyap}, yet to be found. Matano's Ansatz is $L_{pp} = \exp(g(x,u,p))$, which yields a first order PDE for the unknown $g(x,u,p)$ that can be solved by the method of characteristics.
Instead, to accommodate degeneracies, we consider:
\begin{equation}\label{exp}
  L_{pp} := f_q(x,u,p,0,0) \exp({g(x,u,p)}),  
\end{equation}
for some function $g(x,u,p)$ to be found. Note \eqref{exp} is not identically zero, since $f_q(x,u,p,0,0)\not\equiv 0$ due to \eqref{pardeg}.
Moreover, whenever $f_q(x,u,p,0,0)= 0$, the Ansatz \eqref{exp} implies that $L_{pp}\equiv 0$ for any bounded function $g(x,u,p)$. 
However, whenever $g(x,u,p)$ is unbounded, then the interplay between $f_q(x,u,p,0,0)$ and $g(x,u,p)$ plays a major role in the new Ansatz \eqref{exp}, and thus in the construction and regularity of the energies using the present method, in contrast to \cite{Matano88,LappicyFiedler19}. 

To construct the unknown $g(x,u,p)$, we suppose that along the characteristic equations given by
\begin{equation}\label{intro:charac}
\begin{aligned}
    \dot{x}&=f_q(x,u,p,0,0),\\
    \dot{u}&=f_q(x,u,p,0,0) \, p,\\
    \dot{p}&=F^0(x,u,p),
\end{aligned}
\end{equation}
there is a solution $g$ of the following equation:
\begin{equation}\label{intro:Lyapg}
    \dot{g}=-F^0_p(x,u,p) -f_{qx}(x,u,p,0,0)-f_{qu}(x,u,p,0,0) \, p.
\end{equation}
Note that the characteristic equations \eqref{intro:charac} for degenerate PDEs is different from the one obtained by Matano in the non-degenerate case. 
Nevertheless, these equations can be transformed into each other by a suitable `time' rescaling that absorbs $1/f_q(x,u,p,0,0)<\infty$ in case of non-degenerate equations.
Moreover, global existence of the characteristic equations \eqref{intro:charac}-\eqref{intro:Lyapg} might fail, in general. 

\textcolor{black}{Further issues arise when \eqref{diffEQ} is degenerate, i.e. when $f_q(x,u,p,0,0)=0$ for some $(x,u,p)\in [0,1]\times \mathbb{R}^2$, since at such degeneracy points the first two equations of the characteristics \eqref{intro:charac} are $\dot{x}=\dot{u}=0$.
On one hand, if $F^0(x,u,p)\neq 0$, 
then $\dot{p}\neq 0$, which triggers an eventual $f_q(x,u,p,0,0)\neq 0$, due to \eqref{pardeg}. 
On the other hand, if $F^0(x,u,p)=0$, then $\dot{p}=0$ and \eqref{intro:charac} encounters an equilibrium. Hence, the function $g(x,u,p)$ may not be constructed along solutions of \eqref{diffEQ}. 
We call this the \emph{obstacle problem}. 
Therefore, the obstacle problem 
is crucial to rigorously construct energies for degenerate equations using the present method. We explore some of these problems in the examples of Section \ref{sec:ex}.
}
\begin{thm}\label{thm1}
    Assume $f\in C^2$ satisfies \eqref{par}. Suppose the characteristic equations \eqref{intro:charac}--\eqref{intro:Lyapg} have global solutions \textcolor{black}{and that $L_{pp}$ in \eqref{exp} is twice-integrable in $p$}.
    
    Then there exists a Lagrange function $L=L(x,u,p)$ on bounded sets of $(u,p) \in\mathbb{R}^2 $ 
    such that $E:= \int_{0}^{1} L(x,u,u_x) \, dx$
    is a Lyapunov function as in \eqref{IntroLyap} for the equation \eqref{FULLYEQ}. More precisely, bounded solutions $u(t,x)$ of \eqref{FULLYEQ} satisfy
    \begin{equation}
    \label{Lyapdecay}
        \frac{dE}{dt}=-\int_0^1 \exp(g(x,u,u_x)) F^1(x,u,u_x,u_{xx},u_t) \cdot u_t \, dx ,
    \end{equation}
    where $g(\cdot)$ solves \eqref{intro:Lyapg} and \textcolor{black}{$F^1\cdot u_t \geq 0$; the equality holds if, and only if, $u_t\equiv 0$.} 
    
\end{thm}

For non-degenerate quasilinear equations, $    f(x,u,p,q,r)=-r+a(x,u,p)q+h(x,u,p)$, 
where $a>0$, a Lyapunov function $E$ was constructed, independently, by Zelenyak \cite{Zelenyak68} and Matano \cite{Matano88}. 
See also 
\cite{FiedlerRagazzoRocha14} for concise expositions of Matano's method. This method was extended to fully nonlinear non-degenerate parabolic equations, when $f_q\cdot f_r <0$, in \cite{LappicyFiedler19}.
An analogous method for Jacobi systems, a spatially discrete variant, was developed in \cite{FiedlerGedeon99}. For an adaptation to diffusion with singular coefficients see \cite{LappicySing}. 

We emphasize that the procedure to construct the energy function in \eqref{IntroLyap} is formal. Once the Lagrange function $L$ is obtained, one needs to verify various properties needed for a well-defined Lyapunov function, such as integrability, bounds, regularity, etc. 
For this reason, we call the formulae obtained using our method as \emph{energy candidates}.
Thus, the properties and applicability of each candidate still has to be dealt with on a case-by-case basis. 
In Section \ref{sec:ex}, we provide several examples of candidates. 
%
%
%
Note that even if a Lyapunov function is only well-defined for sufficiently regular initial data, one may obtain dynamical information on invariant subspaces of regular enough initial data, \textcolor{black}{see \cite{Souplet,Attouchi,Laurencot,Stinner,ZhangLi}. In particular, non-degenerate equations possess enough regularity to produce a well-defined and regular energy.}
For a deeper regularity analysis of degenerate equations, see \cite{Kalashnikov87,diB,BUV,BDM,TU} and references therein. 
\textcolor{black}{Similarly, energy candidates can potentially be used to obtain local energy estimates akin to \cite{diB}.}



We comment on modifications and possible applications of our result.

Note that our method can potentially treat singular diffusion, i.e., when $f_q(x,u,p,0,0)$ may be unbounded.
An example is $u_t=(u^m u_x/|u_x|)_x$ for $m\geq 0$, which is called the \emph{total variation flow} for $m=0$, or the \emph{heat equation in transparent media} for $m=1$; see \cite{Giacomelli17} and references therein.
However, there are two delicate issues to obtain a Lyapunov function. First, the characteristic equations \eqref{intro:charac} may not have global solutions. 
Second, the Ansatz for $L_{pp}$ defined in \eqref{exp} might not be twice integrable (in $p$) in order to obtain a well-defined formula for $L$ in \eqref{L}. 
\textcolor{black}{A further analysis of singular points must be pursued.}

An alternative splitting of the fully nonlinear equation was pursued in \cite{LappicyFiedler19}, different than \eqref{diffEQ}, yielding an energy that decays according to
\begin{equation}\label{aesthetics}
\frac{dE}{dt}= -\int_{0}^1 L_{pp}\tilde{F}^1  u_t^2 \, dx    
\end{equation}
for some $\tilde{F}^1>0$. 
Instead of the decay in \eqref{Lyapdecay}, one may also be able to obtain a Lyapunov function that decays according to \eqref{aesthetics}, which extracts the $L^2$-gradient flow with weight $L_{pp}\tilde{F}^1>0$.
However, we believe that these different splittings do not change the Lyapunov function itself, only the aesthetics of the abstract formulae.

A semiflow treatment of fully nonlinear degenerate equations \eqref{FULLYEQ} on an appropriate phase-space $X$ 
has been lacking in its full generality, akin to the one for non-degenerate equations provided by \cite{Lunardi95}.
We expect that additional growth conditions on $f$, similar to the non-degenerate case in \cite[Proposition 3.5]{Lunardi91} and \cite[Chapter 6, Sec. 5]{Ladyzhenskaya68}, imply that solutions of \eqref{FULLYEQ} are bounded, global and generate a dissipative semiflow.
In particular, this would guarantee the global existence of the characteristics \eqref{intro:charac}--\eqref{intro:Lyapg} after an appropriate cut-off of $f$ outside a sufficiently large set, and thus the existence of a Lyapunov function $E$ in such a bounded set.
%
In more general settings, including solutions which blow-up, boundedness of $E$ from below may fail. 
In fact, a delicate analysis of the characteristic equations \eqref{quasishootlyap} beyond such crude cut-off may be required in case of blow-up.
\textcolor{black}{
For the non-existence of grow-up (i.e. infinite time blow-up) solutions using such Lyapunov functions, see \cite{Souplet,Attouchi}.}

In addition, it would be desirable to extract dynamic information on the long-term behavior of solutions of \eqref{FULLYEQ}.
Indeed, 
under certain conditions on $f$ that also guarantee asymptotic compactness of the semiflow, there should be a \emph{global attractor} $\mathcal{A}\subset X$ as in the non-degenerate case in \cite{Hale88} or \cite[Theorem 2.2]{Lady91}. 
For particular cases of degenerate type, see \cite{CaChDl,Efendiev}.
%
Thus, as a consequence of the Lyapunov function \eqref{IntroLyap}, bounded trajectories should converge to (sets of) equilibria, according to the LaSalle invariance principle; see \cite[Section 4.3]{Henry81} and \cite[Chapter 5.7]{BabinVishik92} for the non-degenerate case. 
\textcolor{black}{However, the complete description of $\omega$-limit sets is a delicate issue for the degenerate case. See \cite{AronsonCrandall82,FeireislSimonodon99,Carriloetal01,Vazquez04,BonforteFigalli} for specific degenerate cases, which are not in a fully nonlinear setting. In general, see \cite{Hirsch,Matano87,SmithThieme,Polacik,Smith} for a broad overview on the theory of strongly monotone semiflows, when convergence to the set of stationary solutions can be proved.}
Finally, the connection problem for the equations \eqref{FULLYEQ}-\eqref{par} that describes which equilibria are connected by means of a heteroclinic orbit remains open, see \cite{LappicyFully} and references therein for the non-denegerate case.
%
%


The remainder of the paper is organized as follows. In Section \ref{sec:pf}, we prove Theorem \ref{thm1}. In Section \ref{sec:ex}, we compute several significant examples of Lyapunov functions. 

\section{Proof}\label{sec:pf}

We recall the equation \eqref{diffEQ},
\begin{equation}\label{DIFFEQ}
    f_q(x,u,p,0,0)u_{xx} = F^0(x,u,p) + F^1(x,u,p,q,r),
\end{equation}
with  degenerate parabolicity conditions $F^1_r > 0$ and $f_q \geq F^1_q$. See \eqref{par3}.

Differentiating the definition \eqref{IntroLyap} of the Lyapunov function $E$ with respect to time $t$ along classical solutions $u(t,x)$ of \eqref{FULLYEQ}, we obtain
\begin{equation}\label{part1}
    \frac{dE}{dt}= \int_{0}^1 \left( L_u u_t + L_p u_{x t}\right) \, dx .
\end{equation}
Here we used that $u_{xt}=p_t$. 
The Lagrange function $L$ depends on $(x,u,p)=(x,u,u_x)$, only.
It remains to determine $L$ such that $dE/dt<0$, except at equilibria. 
Integrating the term $L_p u_{x t}$ in \eqref{part1} by parts, 
and carrying out the differentiation of $L_p$ with respect to $x$, we obtain
\begin{align}\label{part2}
\begin{split}
    \frac{dE}{dt}&= L_pu_t\Big|_{x=0}^{x=1} +\int_{0}^1 \left( L_u -\frac{d}{dx}L_p  \right) u_t \, dx\\
    &=L_pu_t\Big|_{x=0}^{x=1}+ \int_{0}^1 \left( L_u -L_{px} -L_{pu}u_x -L_{pp}u_{xx}\right) u_t \, dx.
\end{split}
\end{align}
At this point, Matano would plug in the non-degenerate PDE in $u_{xx}$. However, this can not be performed for degenerate equations, since we can not isolate $u_{xx}$ in equation \eqref{DIFFEQ}, as $f_q$ may be zero. In order to remedy this, we modify Matano's original Ansatz, $L_{pp} = \exp(g(x,u,p))$, which would yield a first order PDE to be solved for $g(x,u,p)$. 
Instead, we consider the different Ansatz \eqref{exp} for some function $g(x,u,p)$, yet to be found. 
Note \eqref{exp} is not identically zero, since $f_q(x,u,p,0,0)\not\equiv 0$ due to \eqref{par}.
Thus
\begin{align}\label{part2b}
        \frac{dE}{dt}&=L_pu_t\Big|_0^1 +\int_{0}^1 \left( L_u -L_{px} -L_{pu}u_x - \exp({g}) f_q u_{xx}\right) u_t \, dx.
\end{align}
We then substitute the PDE \eqref{FULLYEQ} recast in \eqref{DIFFEQ}, to obtain
\begin{equation}\label{part2c}
        \frac{dE}{dt}=L_pu_t\Big|_0^1 +\int_0^1(L_u -L_{px} -L_{pu}u_x - \exp({g}) F^0)u_t \, dx - \int_0^1 \exp({g}) F^1u_t\, dx.
\end{equation}

We seek to construct the Lagrange function $L$ such that the boundary terms vanish, the parenthesis in the first integral \eqref{part2c} also vanishes, and satisfies the Ansatz \eqref{exp} for some function $g(x,u,p)$. This yields a Lyapunov function such that
\begin{equation}\label{Lyaput2}
        \frac{dE}{dt}=-\int_{0}^1 \exp({g}) F^1 u_t \, dx.
\end{equation}
Note $F^1u_t\geq 0$, due the parabolicity condition $F^1_r>0$.
Next, we guarantee that there exists a function $g(x,u,p)$ such that 
\begin{equation}\label{Lyapquasi0}
    L_u-L_{px}-pL_{pu} -\exp({g}) F^0=0,
\end{equation}
for all $(x,u,p)\in [0,1]\times \mathbb{R}^2$, and also $L_pu_t=0$ on the boundaries $x=0,1$.
Note that $(u,p)\in \mathbb{R}^2$ are real variables rather than solutions $u,u_x$ of PDEs depending on $(t,x)$.
Differentiating \eqref{Lyapquasi0} with respect to $p$, the terms $L_{pu}$ cancel, yielding
\begin{equation}\label{Lyapquasi0_p}
    L_{ppx}+pL_{ppu}+\exp({g}) g_p F^0 =-\exp(g)F^0_p.
\end{equation}
Rewriting \eqref{Lyapquasi0_p} in terms of $g$, according to \eqref{exp}, amounts to the first order PDE,
\begin{equation}\label{Lyapquasig}
    f_qg_{x}+pf_qg_{u}+F^0g_{p} =-F^0_p -f_{qx}-pf_{qu} .
\end{equation}
The method of characteristics can solve \eqref{Lyapquasig}: along solutions of the auxiliary ODEs
\begin{equation}
\label{quasishootlyap}
\begin{aligned}
    \dot{x}&=\frac{dx}{d\tau}=f_q(x,u,p,0,0),\\
    \dot{u}&=\frac{du}{d\tau}=f_q(x,u,p,0,0) \, p,\\
    \dot{p}&=\frac{dp}{d\tau}=F^0(x,u,p),
\end{aligned}
\end{equation}
the function $g$ must satisfy
\begin{equation}\label{Lyapg}
    \dot{g}=\frac{dg}{d\tau}=-F^0_p(x,u,p) -f_{qx}(x,u,p,0,0)-f_{qu}(x,u,p,0,0) \, p,
\end{equation}
with the initial condition $g(0,u_0,p_0)$, where $(u_0,p_0):=(u(0,0),u_x(0,0))$. Our differentiability assumptions on $f$ imply $g \in C^0$, at least. 
 
Without further assumptions on the nonlinearity $f$ in \eqref{FULLYEQ}, solutions to \eqref{quasishootlyap} may not exist on the whole required interval $x\in [0,1]$. 
For this reason, we have assumed the global existence of solutions for the characteristic equations.
\textcolor{black}{Moreover, note that the global existence of the characteristics is not enough to guarantee the existence of a Lyapunov function, in general.
Indeed, further complications occur when the diffusion degenerates (i.e. $f_q(x,u,p,0,0)=0$), since $\dot{x}=\dot{u}=0$ in \eqref{quasishootlyap}.
If $F^0(x,u,p)\neq 0$, then $\dot{p}\neq 0$ and this may trigger an eventual $f_q(x,u,p,0,0)\neq 0$. 
However, if $F^0(x,u,p)=0$, then \eqref{quasishootlyap} encounters an equilibrium and the function $g(x,u,p)$ can not be constructed along solutions of \eqref{quasishootlyap}. 
We call this the \emph{obstacle problem}. 
\textcolor{black}{Therefore, the obstacle problem 
is crucial to rigorously construct energies for degenerate equations using the present method. We explore some of these issues in the examples of Section \ref{sec:ex}.}
}

After this construction, we now have to reverse gear and ascend from a function $g$ satisfying \eqref{Lyapquasig} to a Lagrange function $L$ satisfying \eqref{Lyapquasi0}. The general solution $L$ of $L_{pp}=f_q\exp(g)$ can be obtained by integrating it twice with respect to $p$:
\begin{equation}\label{L}
    \begin{aligned}
        L(x,u,p):=&\int_0^p \int_0^{p_1} f_q(x,u,p_2,0,0)\exp(g(x,u,p_2))\, dp_2\, dp_1 \\
        &+ L^0(x,u)+L^1(x,u)p.
    \end{aligned}
\end{equation}
This solves \eqref{Lyapquasi0_p}. 
To ensure that $L$ is also a solution of \eqref{Lyapquasi0}, we have to determine the integration ``constants'' $L^0$ and $L^1$, appropriately. 
%
Recall that \eqref{Lyapquasi0_p} was obtained through differentiation of \eqref{Lyapquasi0} with respect to $p$. Conversely, the left-hand side of \eqref{Lyapquasi0} is therefore independent of $p$. Hence \eqref{Lyapquasi0} is satisfied for all $p$, if it holds for some fixed value $p=p_*\in\mathbb{R}$. 

{\color{black} 
Deriving \eqref{L} with respect to $u$ and $p$ yields
\begin{subequations}
    \begin{align}
        L_u &= \int_0^{p} \int_0^{p_1} \left( f_{qu}
        + f_q
        g_u
        \right)\exp(g
        )\, dp_2\, dp_1 + L_u^0(x,u)+L_u^1(x,u)p,\label{Lu}\\
        L_p &= \int_0^{p} f_q
        \exp(g
        ) dp_1+ L^1(x,u),\label{Lp}
    \end{align}
\end{subequations}
where the integrand arguments 
are suppressed to alleviate the notation. Moreover, further differentiating $L_p$ with respect to $x$ and $u$ produces
\begin{subequations}\label{LpxLpu}
    \begin{align}
        L_{px}&=\int_0^{p} \left( f_{qx} +f_{q} g_x \right) \exp(g) dp_1+ L^1_x(x,u),\\
        L_{pu}&=\int_0^{p} \left(  f_{qu} +f_{q} g_u \right) \exp(g) dp_1+ L^1_u(x,u).
    \end{align}
\end{subequations}

Note that evaluating the equation \eqref{Lyapquasi0} at $p_
*$ yields $L_u=L_{px}+p_*L_{pu}+\exp({g}) F^0$. Substituting \eqref{Lu} and \eqref{LpxLpu}, evaluated at $p_*$, and isolating $L^0_u$, we obtain that
\begin{equation}\label{Lu^0}
    \begin{aligned}
        L_u^0(x,u)  = & L^1_x(x,u) 
        + \exp(g(x,u,p_*))F^0(x,u,p_*) \\ 
        &+\int_0^{p_*} \left[ f_{qx} +f_{q} g_x  + \left( f_{qu} +f_{q} g_u\right)p_* \right]\exp(g) dp_1\\
        &- \int_0^{p_*}\int_0^{p_1} \left( f_{qu}+ f_q g_u \right)\exp(g)\, dp_2 dp_1.  
    \end{aligned}
\end{equation}
%
Note that the right hand side of \eqref{Lu^0} depends only on $f,g,L^1$ evaluated at $p_*$.
Hence, we can integrate \eqref{Lu^0} with respect to $u$, which in turn yields the function $L^0(x,u)$ explicitly written as a function of $(x,u)$ for any $p_*$. Mathematically, we achieve that
\begin{equation}\label{L^0}
    \begin{aligned}
        L^0(x,u)  = \int_0^u \Bigg[ & L^1_x(x,u_1) + \exp(g(x,u_1,p_*))F^0(x,u_1,p_*) \\ 
        &+\int_0^{p_*} \left[ f_{qx} +f_{q} g_x  + \left( f_{qu} +f_{q} g_u\right)p_* \right]\exp(g) dp_1\\
        &- \int_0^{p_*}\int_0^{p_1} \left( f_{qu}+ f_q g_u \right)\exp(g)\, dp_2 dp_1 \Bigg] du_1 + L^{00}(x).  
    \end{aligned}
\end{equation}
%
%
To complete the proof, it only remains to show that $L_pu_t$ vanishes at the boundaries $x=0,1$, which is done by appropriately constructing $L^1$. 
At any boundary of Dirichlet type \eqref{Dir} this is trivial because $r=u_t=0$. Thus we can either let $L^1\equiv 0$, or
\begin{equation}\label{L1dir}
    L^1(x,u) :=  -  \int_0^{p_*} f_q(x,u,p,0,0)\exp(g(x,u,p))dp,
\end{equation}
which respectively yields that $L_p$ is a finite value or zero, according to \eqref{Lp}. Note that the choice of $L^1$ influences the construction of $L^0$ in \eqref{L^0}. 

In the case of a nonlinear Robin boundary condition \eqref{Rob} at only one boundary, either $x=0$ or $x=1$, we have to choose $L$ such that $L_p(\iota,u,b^\iota(u))=0$. 
By our construction \eqref{L} of $L$, on behalf of \eqref{Lp}, this is equivalent to 
\begin{equation}\label{G1}
    L^1(\iota,u):=-\int_0^{b^\iota(u)} f_q(\iota,u,p,0,0) \exp(g(\iota,u,p))\, dp,
\end{equation}
and we may choose $L^1$ to be independent of $x$.

For nonlinear Robin boundary conditions \eqref{Rob} at both boundaries, $x=0$ and $x=1$, we define $L^1(\iota,u)$ as in \eqref{G1} for $\iota=0,1$. 
Therefore, the linear interpolation $L^1(x,u):=(1-x)L^1(0,u)+xL^1(1,u)$ provides $L^1 \in C^1$ such that $L_p(\iota,u,b^\iota(u))=0$.
} 

For example, if $p_*=0$, the construction of $L$ yields $L_p=L^1$, $L_{px}=L^1_x$ and $L_u=L^0_u$. Evaluating either \eqref{Lyapquasi0} or \eqref{Lu^0} at $p_*=0$ yields $L^0_u=L^1_x+\exp(g) F^0$.
Integrating with respect to $u$, \textcolor{black}{in agreement with \eqref{L^0}, } we can neglect an irrelevant additive constant $L^{00}(x)$ for $E$ to obtain that
\begin{equation}\label{G0}
    L^0(x,u):=\int_0^{u} \left[L^1_x(x,u_1) 
    +\left( \exp(g(x,u_1,p_*))F^0(x,u_1,p_*)\right)|_{p_*=0}\right]du_1.
\end{equation} 

\textcolor{black}{For $p_*=0$, the choices of $L^1$, which depend on the boundary conditions, yield $L^1\equiv 0$ for Dirichlet boundary conditions and \eqref{G1} is unchanged for Robin boundary conditions. }




\section{Examples}\label{sec:ex}

We explicitly compute examples of energy candidates using the method in the previous section and compare them with well-known Lyapunov functions in the literature; see Table \ref{tab}. 
For the sake of simplicity, we consider Dirichlet boundary conditions throughout the examples, which yield $L^1\equiv 0$.
%
\begin{table}[H]
    \color{black}
    \scriptsize
    \centering
    \begin{tabular}{|c|c|l|}
        \rowcolor{Gray} \hline
        Section & PDE & \hspace{4.6cm} Energy \\ \hline
        \ref{sec:gradDIFFhamiltREAC} & $u_t=a(u_x)u_{xx}+h(u)$ & 
        \begin{tabular}{cc}
            \underline{Old}: & $\displaystyle{\int_0^1}\left( \int_0^p \int_0^{p_1} a(p_2) d{p}_2 d{p}_1 -\int_0^u h(u_1) du_1\right)dx$ \vspace{0.3cm} \\
            \underline{New}: & Same
        \end{tabular}
        \\ \hline
        \ref{sec:ICMF} & $u_t = \frac{1+u^2_{x}}{1 -\left(1- \frac{u_x^2}{1+u_x^2}\right)u_{xx}}$ & 
        \begin{tabular}{cl}
            \underline{Old}: & $|M_t^2|^{1/2}\left(8 \pi - \displaystyle{\int_{M_t^2}}H^2 d\mu_t\right)$, see \cite{Marquardt13},
              \cite[Prop. 6.1]{Marquardt17} \vspace{0.3cm} \\
            \underline{New}: & $\displaystyle{\int_0^1} u_x\arctan (u_x) - \log(1+u_x^2) -  u  \, dx$
        \end{tabular}
        \\ \hline
        \ref{sec:MCFreac} & $ u_t = \left(\frac{u_{x}}{\sqrt{1 + u_x^2}}\right)_x + u_x^n$ & 
        \begin{tabular}{cc}
            \underline{Old}: & Unknown \vspace{0.3cm} \\
            \underline{New}: & \begin{tabular}{ll}
                    &$\displaystyle{\int_{0}^{1}} \left(\displaystyle{\int_0^p}\displaystyle{\int_0^{p_1}}\dfrac{1}{|p_2|^n(1 + p_2^2)^{\frac{3}{2}}} dp_2 \ dp_1 \right) - u \, dx$\\ 
                    & See Table \ref{tab2} for $n=1,2,3,4,5$
                    \end{tabular}
        \end{tabular}
        \\ \hline
        \ref{sec:rhopoly} & $u_t = (|u_x|^{\rho-2} u_x)_x+u_x^n$ & 
        \begin{tabular}{cc}
            \underline{Old}: & \begin{tabular}{cc}
                    \underline{$\rho=2$}: & $\displaystyle{\int_0^1} \bigg(\dfrac{|u_x|^{2 - n}}{(2-n)(1-n)} - u \bigg)\, dx$, see \cite{Souplet,Laurencot}\\
                    \underline{$\rho\neq 2$}: & Unknown, see \cite{Attouchi,Stinner}
                   \end{tabular} \vspace{0.3cm} \\
            \underline{New}: & \begin{tabular}{lc}
                    \underline{$n\neq \rho,\rho-1$}: & $\displaystyle{\int_0^1} \ \dfrac{(\rho - 1)}{(\rho -n)(\rho - n-1)} |u_x|^{\rho - n} - u\ dx$  \\
                    \underline{$n=\rho-1$}: & $\displaystyle{\int_0^1} (\rho-1)|u_x| \left(\log|u_x|-1\right) -u\ dx$\\
                    \underline{$n = \rho$}: & $\displaystyle{\int_0^1} (1 - \rho)\log|u_x|-u \ dx$
                    \end{tabular}
        \end{tabular}
        \\ \hline        
        \ref{sec:PME} & $u_t =(u^m)_{xx}$ & 
        \begin{tabular}{cl}
            \underline{Old}: & $\displaystyle{\int_0^1} \frac{u^{m+1}}{m+1}dx$, \text{see} \cite{AronsonCrandall82,DiazDiaz,Vazquez06},\cite[Eq. (2.7)]{Vazquez04} \vspace{0.3cm} \\
            \underline{New}: & $\displaystyle{\int_0^1} m u^{m-1} |u_x|\left(\log |u_x| - 1\right) dx$
        \end{tabular}
         \\ \hline
    \end{tabular}
    \caption{\textcolor{black}{Comparison of known Lyapunov functions and the new candidate formulae for specific PDEs using Matano's method.}}
    \label{tab}
\end{table}

\subsection{Gradient-degenerate quasilinear diffusion with nonlinear forcing}\label{sec:gradDIFFhamiltREAC}
Consider the equation
\begin{equation}\label{degPDE1}
    u_t = a(u_x)u_{xx}+h(u), 
\end{equation}
with \textcolor{black}{$a,h\in C^2$ such that $a(u_x) \geq 0$, where the equality only happens in a set of measure zero, due to \eqref{pardeg}}. 
In the abstract setting in the previous section, we have that $f_q= a(p)$, $F^0 = -h(u)$ and $F^1 = u_t$. 

Thus, the characteristic equations \eqref{quasishootlyap} are given by
\begin{equation}\label{degshoot}
    \begin{aligned}
    \dot{x}&= a(p),\\
    \dot{u}&= a(p)\, p,\\
    \dot{p}&= -h(u),
\end{aligned}
\end{equation}
and $g$ evolves according to \eqref{Lyapg}, i.e.,
\begin{equation}
    \dot{g}= 0.
\end{equation}
{\color{black}
Note that if $h(u)\neq 0$, then 
$\dot{p}\neq 0$ and the equation \eqref{degshoot} does not encounter an equilibrium obstacle. 
However, if $h(u)=0$ for some $u\in\mathbb{R}$, then there is a constant equilibrium of the PDE \eqref{degPDE1}, which is also an equilibrium obstacle for the characteristic equations. 
However, in either case, note that $\dot{g}=0$. 
}
Therefore $g(x(\tau),u(\tau),p(\tau))$ is a constant function along any solution of the \eqref{degPDE1}, and therefore we obtain the trivial solution $g \equiv 0$ for the initial condition $g_0:=g(0,u(0),p(0)) = g(0,u_0,p_0) = 0$. 
Due to the equation \eqref{exp}, we obtain that $L_{pp} = a(p)$.
Note that $L^1\equiv 0$ due to Dirichlet boundary conditions
, and $L^0=-\int_0^u h(u_1)\, du_1$ due to the equation \eqref{G0} for $p_*=p_0$. Hence \eqref{L} implies that
\begin{equation}\label{Lyap1}
E = \int_0^1\bigg( \int_0^p \int_0^{p_1} a(p_2)\ d{p}_2 \ d{p}_1 -\int_0^u h(u_1)\, du_1\bigg) \, dx,
\end{equation}
which decays according to 
\begin{equation}\label{decPLAP}
    \frac{dE}{dt} = -\int_0^1 u_t^2 \, dx.
\end{equation}
\textcolor{black}{Note that for Robin boundary conditions, equation \eqref{G1} yields the following term, $L^1(x,u) = -(1-x)\int_0^{b^0(u)} a(p_1)dp_1 - x\int_0^{b^1(u)}a(p_1)dp_1$, whereas equation \eqref{G0} implies $L^0(x,u) = \int_0^u \left(\int_{b^1(u_1)}^{b^0(u_1)} a(p_1)dp_1 - h(u_1)\right) du_1$. Therefore, the energy candidate formula \eqref{Lyap1} can be modified accordingly. We reiterate that the rigorousness of this formula depends on the delicacy of solutions of the characteristic equations \eqref{degshoot}. }

In particular, the $\rho$-Laplacian\footnote{In the literature, this operator is called the $p$-Laplacian. However, in our notation $p:=u_x$ and thus we replace the parameter $p$ by $\rho$ in the degenerate diffusion operator, i.e., $\partial_\rho u:=(|u_x|^{\rho-2}u_x)_x$.} equation occurs when $a(u_x) = (\rho-1)|u_x|^{\rho - 2}$ and $h\equiv 0$, and thus we recover its well-known energy $E = \int_0^1 |u_x|^{\rho}/\rho \, dx$.
Also, the mean curvature flow for one dimensional graphs occurs when $a(u_x) = (1 + u_x^2)^{-3/2}$ and $h\equiv 0$, and thereby we also recover the energy $E = \int_0^1 \sqrt{1+u_x^2} \, dx$. This energy accounts for the perimeter of the curve, which decreases under evolution of mean curvature according to \eqref{decPLAP}. 
\textcolor{black}{For a proof of infinte time blow-up (i.e. grow-up) for a mean curvature flow with a general Hamiltonian reaction of type $h(x,u)$, due to the existence of this well-known energy formula, see \cite{Chen}. }
Note that the mean curvature flow only degenerates at infinity, i.e., when $|u_x|\to \infty$, and thus we expect that an appropriate compactification of the semiflow will be described by a degenerate equation at infinity, 
see \cite{LappicyJu}.

\textcolor{black}{The energy formula \eqref{Lyap1} is ubiquitous in the theory of degenerate parabolic PDEs in divergence form, whenever the degeneracy occurs in the gradient. It is a well-defined and sufficiently regular energy with several consequences which include well-posedness, regularity and dynamical properties. See, for example, \cite{Efendiev, CaChDl, diB}}.
%
%
%
\subsection{Inverse mean curvature flow for certain graphs}\label{sec:ICMF}
Consider the equation
\begin{equation}\label{IMCF}
    u_t = \frac{1+u^2_{x}}{1 -\left(1- \frac{u_x^2}{1+u_x^2}\right)u_{xx}}.
\end{equation}
This equation has been considered in higher dimensions in \cite[Section 3]{Marquardt13}, and we construct a different monotone quantity in comparison to \cite[Proposition 6.1]{Marquardt17}.

In this case, we have that
\begin{equation}\label{IMCFF1}
    f_q = \frac{(1 + u_x^2)^2}{(1 + u_x^2-u_{xx})^{2}}, \quad F^0=-(1+u_x^2), \quad F^1=u_t +  \frac{u^2_{xx}}{u_{xx} - (1+u_x^2)}
\end{equation}
%
%
Note that this equation is not degenerate, since $f_q(x,u,p,0,0)=1$, \textcolor{black}{but it is singular whenever $u_{xx}=u_x^2+1$}.
The characteristic equations\eqref{quasishootlyap} is given by
\begin{equation}
    \begin{aligned}
    \dot{x} &= 1,\\
    \dot{u} &= p,\\
    \dot{p} &= -(1+p^2),
    \end{aligned}\label{characIMCF}
\end{equation}
and \eqref{Lyapg} is given by 
\begin{equation}
    \dot{g}= 2p.
\end{equation}
\textcolor{black}{Since the equation \eqref{IMCF} is non-degenerate, then $\dot{x}>0$ and the characteristic equations \eqref{characIMCF} does not encounter an equilibrium obstacle. Thus, the global existence of characteristics is enough to pursue the construction in the previous section and guarantee the existence of a Lyapunov function.}
We can solve these equations explicitly:
\begin{equation}\label{pIMCF}
    p(\tau)=-\tan (\tau+\arctan(-p_0))
\end{equation}
and
\begin{equation}\label{gIMCF}
    g(\tau)=g_0+2\log\left(\frac{\cos(\tau+\arctan(-p_0))}{\cos(\arctan(-p_0))}\right).
\end{equation}
Consequently,
\begin{equation}
    \begin{aligned}
        g(p)&=g_0+ 2\log\left(\frac{\cos(\arctan(-p))}{\cos(\arctan(-p_0))} \right),\\
        &=g_0+\log\left( \frac{1+p_0^2 }{1+p^2} \right)
    \end{aligned}
\end{equation}
\textcolor{black}{Note that we obtain a finite value $g(p)$ for any $p\in\mathbb{R}$, including $p=0$. 
}
Thus equations \eqref{exp}, \eqref{G0} with $p_*= 0$, and the Dirichlet boundary imply that 
\begin{equation}
\begin{split}
  L_{pp} = \exp(g_0)\frac{1+p_0^2}{1+p^2} 
  \quad L^0 = - \exp(g_0) (1+p_0^2) u \quad \text{and} \quad L^1 = 0.
    \end{split}
\end{equation}
Thus the energy \eqref{L}, up to a multiplicative constant $\exp(g_0)(1+p_0^2)$, is given by
\begin{equation}\label{E:IMCF}
    E= \bigintssss_{0}^{1} \bigg(
    u_x\arctan (u_x) - \log(1+u_x^2) -  u \bigg) \, dx,
\end{equation}
which decays according to
\begin{equation}\label{E:IMCFdecay}
    \frac{dE}{dt} = - \bigintssss_{0}^{1} \frac{(2+u_x^2)u_x^2}{(1+u_x^2)^3}u_t^2 \, dx,
\end{equation}
since $F^1$ given by \eqref{IMCFF1} can be rewritten as $F^1=\frac{(2+u_x^2)u_x^2}{(1+u_x^2)^2}u_t$ by substituting \eqref{IMCF}.
\textcolor{black}{Note that neither the energy formula \eqref{E:IMCF}, nor its decay in \eqref{E:IMCFdecay}, possess singularities for bounded values $u,u_x\in\mathbb{R}$. }

\textcolor{black}{Note that for Robin boundary conditions, equation \eqref{G1} yields the following term, $L^1(x,u) = -\exp(g_0)(1+p_0^2)\left[(1-x)\arctan(b^0(u)) + x \arctan(b^1(u))\right]$, whereas equation \eqref{G0} implies $L^0(x,u) = \exp (g_0)(1+p_0^2)\int_0^u \arctan(b^0(u_1)) - \arctan(b^1(u_1)) \ du_1$. Therefore the energy formula \eqref{E:IMCF} can be modified accordingly. }

We emphasize that, in this example, we have computed an energy formula for a fully nonlinear non-degenerate equation, where the method in \cite{LappicyFiedler19} is not applicable. 
Indeed, the splitting of the PDE \eqref{IMCF} according to \cite{LappicyFiedler19} defines different functions $F,F^0,F^1$ than our present method, which yields a different splitting of the PDE in contrast to \eqref{diffEQ}. 
In particular, the functions $F,F^0,F^1$ in \cite{LappicyFiedler19} are not well-defined for this example.
Roughly speaking, trying to isolate $u_{xx}$ in \eqref{IMCF} to define the function $F$ in \cite{LappicyFiedler19} yields an ill-defined vector field when $u_t=0$.
Therefore, the present example shows that our current method, which splits the PDE according to \eqref{diffEQ}, overcomes certain problems arising even in the non-degenerate construction in \cite{LappicyFiedler19}. 

\subsection{Mean curvature flow with an external forcing}\label{sec:MCFreac}
Consider the equation \textcolor{black}{that describes the mean curvature flow for planar graphs with an external forcing given by $u_x^n$ with $n\in \mathbb{N}$,}
\begin{equation}\label{MCFpoly}
    u_t = \left(\frac{u_{x}}{\sqrt{1 + u_x^2}}\right)_x + u_x^n= \frac{u_{xx}}{(1 + u_x^2)^{\frac{3}{2}}} + u_x^n.
\end{equation}
In this case, we have that $f_q = (1 + u_x^2)^{-3/2}$, $F^0 = -u_x^n$ and $F^1 = u_t$. \textcolor{black}{Note that for solutions which do not blow-up in the gradient (e.g. differentiable solutions), $f_q>0$ and thereby the equation \eqref{MCFpoly} is non-degenerate. } 
Hence the characteristic equations \eqref{quasishootlyap} are given by
\begin{equation}
    \begin{aligned}
    \dot{x} &= \frac{1}{(1 + p^2)^{\frac{3}{2}}},\\
    \dot{u} &= \frac{p}{(1 + p^2)^{\frac{3}{2}}},\\
    \dot{p} &= -p^n,
    \end{aligned}\label{characMCF}
\end{equation}
and \eqref{Lyapg} is given by 
\begin{equation}\label{gMCF}
    \dot{g}= np^{n-1}.
\end{equation}
%
%
\textcolor{black}{Since the equation \eqref{MCFpoly} is non-degenerate, then $\dot{x}>0$ and the characteristic equations \eqref{characMCF} does not encounter an equilibrium obstacle for finite $p\in\mathbb{R}$. 
Note that if $p_0=0$, then $p(\tau)\equiv 0$ and $g(\tau)\equiv g_0$.\footnote{\textcolor{black}{For $p_0= 0$ (i.e. $p\equiv 0$ and $g\equiv g_0$), we obtain that $L_{pp} = \exp(g_0)/(1 + p^2)^{\frac{3}{2}}$ for all $n\in\mathbb{N}$. 
This yields $E = \int_0^1 \sqrt{1+u_x^2} \, dx$, up to a multiplicative constant $\exp(g_0)$, which is the perimeter of the curve; similar to the mean curvature flow with Hamiltonian forcing in Section \ref{sec:gradDIFFhamiltREAC}.} }
Moreover, if $p_0>0$, then $p(\tau)$ decreases to $0$ as $\tau\to \infty$, however, if $p_0<0$, then $p(\tau)$ either decreases to $0$ as $\tau\to \infty$ or blows up in finite time, respectively for $n$ odd or even.
\textcolor{black}{In addition, if $p(\tau)$ blows up in finite time, note that the characteristic equations encounters a \emph{singularity obstacle}, where $\dot{x}=\dot{u}=0$, but $\dot{p}=\pm\infty$.}
Thus, for $n$ odd, the global existence of characteristics is enough to pursue the construction in the previous section and guarantee the existence of a Lyapunov function, which is not the case for $n$ even.
}


We can solve these equations explicitly,
\begin{equation}\label{pLapPolysol}
    p(\tau)=
    \begin{cases}
        p_0e^{-\tau}\, &\text{ for } n=1\\
        \frac{p_0}{\left(1+(n-1)p_0^{n-1}\tau\right)^{\frac{1}{n-1}}} \, & \text{ for } n>1,
    \end{cases}
\end{equation}
and
\begin{equation}\label{gpoly}
    g(\tau)=
    \begin{cases}
        \tau + g_0\, &\text{ for } n=1\\
         \log\left(1 +(n-1)p_0^{n-1} \tau\right)^\frac{n}{n-1}+g_0\, & \text{ for } n>1.
    \end{cases}
\end{equation}
Consequently,
\begin{equation}
    g(p) = \log\left(\frac{p_0}{p}\right)^n+g_0
\end{equation}
for all $n\in\mathbb{N}$.
%
Therefore, the equations \eqref{exp}, \eqref{G0} with $p_*=p_0\neq 0$, and Dirichlet boundary conditions yield 
\begin{equation}\label{LppMCF}
\begin{split}
  L_{pp} = \exp(g_0)|p_0|^n\dfrac{1}{|p|^n(1 + p^2)^{\frac{3}{2}}} 
    \quad L^0 = - \exp(g_0) |p_0|^n u\quad \text{and} \quad L^1 = 0.
    \end{split}
\end{equation}
Thus the energy \eqref{L}, up to a multiplicative constant $\exp(g_0)|p_0|^n$, is formally given by
{\color{black}
\begin{equation}\label{E:MCFpoly}
    E= \displaystyle{\int_{0}^{1}} \left(\displaystyle{\int_0^p}\displaystyle{\int_0^{p_1}}\dfrac{1}{|p_2|^n(1 + p_2^2)^{\frac{3}{2}}} dp_2 \ dp_1 \right) - u \, dx,
\end{equation}
which decays according to
\begin{equation}
    \frac{dE}{dt} = - \bigintssss_{0}^{1} \frac{u_t^2}{|u_x|^n} \, dx.
\end{equation}
However, note that \eqref{LppMCF} may not be twice integrable for all $n\in\mathbb{N}$ and therefore the energy \eqref{E:MCFpoly} may be ill-defined for some $n\in\mathbb{N}$. This is in contrast with the example in Section \ref{sec:rhopoly}, which possess the same characteristic equation, but it has a different Ansatz for $L_{pp}$.

Since equation \eqref{MCFpoly} is non-degenerate for bounded solutions, one can compute a Lyapunov function for classical bounded solutions following the original construction of Matano in \cite{Matano88,LappicyFiedler19}. 
We now proceed with the original construction to compare with our present results. 
Indeed, instead of separating the equation \eqref{MCFpoly} according to \eqref{diffEQ}, which amounts to the characteristic equations \eqref{intro:charac}, the characteristics in \cite{Matano88,LappicyFiedler19} are given by
\begin{equation}
    \begin{aligned}
    \dot{x} &= 1,\\
    \dot{u} &= p,\\
    \dot{p} &= -p^n (1 + p^2)^{\frac{3}{2}},
    \end{aligned}\label{characMCFv2}
\end{equation}
whereas the unknown $g$ should satisfy the follow equation, in contrast to \eqref{intro:Lyapg}, 
\begin{equation}
    \dot{g}= p^{n-1} [n+(n+3)p^2] \sqrt{1 + p^2}.
\end{equation}
Note that the evolution equation for $p$ in \eqref{characMCFv2} decouples from $(x,u)$ similar to \eqref{characMCF}. Moreover, the global existence of the characteristics equation \eqref{characMCFv2} still depends on the parity of $n$, i.e., global existence only occurs for odd $n$. However, a seemingly more complicated vector field appears in the right hand side. In order to obtain a simpler equation for $p$, which can be solved explicitly, we introduce a new time variable, $\tilde{\tau}$, such that $p':=dp/d\tilde{\tau}=(1+p^2)^{-3/2}\dot{p}$, which transforms the equation \eqref{characMCFv2} into\footnote{\textcolor{black}{Note the characteristics \eqref{characMCFv3} coincide with the one obtained through our construction, see \eqref{characMCF}. However, the equations for $g$ given by \eqref{gMCF3} and \eqref{gMCF} are different. This occurs since our Ansatz in \eqref{exp} is different than Matano's, which is $L_{pp}=\exp(g)$.}}
\begin{equation}
    \begin{aligned}
    x' &= \frac{1}{(1 + p^2)^{\frac{3}{2}}},\\
    u' &= \frac{p}{(1 + p^2)^{\frac{3}{2}}},\\
    p' &= -p^n ,
    \end{aligned}\label{characMCFv3}
\end{equation}
whereas the unknown $g$ satisfies 
\begin{equation}\label{gMCF3}
    g'= \frac{p^{n-1} [n+(n+3)p^2]}{1 + p^2}.
\end{equation}
Similarly to \eqref{pLapPoly}, we can solve the relevant part of the equations \eqref{characMCFv2} explicitly: 
\begin{equation}
    p(\tilde{\tau})=
    \begin{cases}
        p_0e^{-\tilde{\tau}}\, &\text{ for } n=1\\
        \frac{p_0}{\left(1+(n-1)p_0^{n-1}\tilde{\tau}\right)^{\frac{1}{n-1}}} \, & \text{ for } n>1,
    \end{cases}
\end{equation}
and 
\begin{equation}
    g(\tilde{\tau})=
    \begin{cases}
        g_0 + \log (1+p_0^2)^{\frac{3}{2}} + 4\tilde{\tau} + \log\left(\frac{e^{-3\tilde{\tau}}}{(1 + p_0^2e^{-2\tilde{\tau}})^{\frac{3}{2}}}\right)   \, &\text{ for } n=1\\
        g_0 + \log (1+p_0^2)^{\frac{3}{2}} + \log\left(\frac{(1 + \tilde{\tau}(n-1)p_0^{n-1})^{\frac{n+3}{n-1}}}{\left(p_0^2 + (1 + \tilde{\tau}(n-1)p_0^{n-1})^{\frac{2}{n-1}}\right)^{\frac{3}{2}}}\right) \, & \text{ for } n>1.\\
    \end{cases}
\end{equation}
Consequently, 
\begin{equation}
g(p) =  g_0 + \log (1+p_0^2)^{\frac{3}{2}} + \log\left(\frac{p_0^{n}}{p^n(1+p^2)^{\frac{3}{2}}}\right).
\end{equation}
In turn, Matano's Ansatz yields
\begin{equation}\label{MCFMatano}
    L_{pp} = \exp(g) =  \exp(g_0)|p_0|^n(1+p_0^2)^{\frac{3}{2}}\frac{1}{|p|^n(1 + p^2)^{\frac{3}{2}}}.
\end{equation}
Therefore, we compare the resulting $L_{pp}$ in \eqref{MCFMatano} following Matano's construction to the resulting $L_{pp}$ in \eqref{LppMCF} using our construction. Note these are equal, up to a multiplicative constant $(1+p_0)^{3/2}$, and hence both our methods yield the same energy formulae.
Therefore, the lack of integrability of $L_{pp}$ for some $n$ occurs both in our and Matano's methods.
%
In particular, we have used the software \emph{Maple} to formally compute the integrals \eqref{E:MCFpoly} for $n=1,2,3,4,5$, which is displayed in the Table \ref{tab2}.

\begin{table}[H]
    \footnotesize
    \color{black}
    \centering
    \begin{tabular}{|c|l|}\hline
        \rowcolor{Gray} $n$ & \hspace{4.6cm} Energy formula
        \\ \hline 
         $1$& $ E= \displaystyle{\int_0^1} - u_x \tanh^{-1}\left(\sqrt{1 + u_x^2}\right) - u\ dx$\\ \hline
             $2$& $E=\displaystyle{\int_0^1}\tanh^{-1}\left(\sqrt{1+u_x^2}\right) -2\sqrt{1+u_x^2} - u\ dx$\\
         \hline
        $3$& $E=\displaystyle{\int_0^1} \frac{1}{|u_x|}\left(\sqrt{1 + u_x^2} + 3u_x^2 \tanh^{-1}{\left(\sqrt{1 + u_x^2}\right)} \right)\ - u \ dx $\\
         \hline
             $4$ & $E=\displaystyle{\int_0^1} \frac{(1+16u_x^2)\sqrt{1+u_x^2}}{ u_x^2} - \frac{\tanh^{-1}\left(\sqrt{1+u_x^2}\right)}{2} - \frac{u}{3^{1/3}} \ dx.$\\
         \hline
         $5$& $E=\displaystyle{\int_0^1}-\left(\frac{1}{u_x^2|u_x|}\left(\sqrt{1 + u_x^2}(-2 + 19u_x^2) + 45u_x^4\tanh^{-1}(\sqrt{1+u_x^2})\right) + u\right) \ dx $\\
         \hline
    \end{tabular}
    \caption{\textcolor{black}{Explicit examples of the energy candidates \eqref{E:MCFpoly} for the equation \eqref{MCFpoly} for $n=1,2,3,4,5$. Note that for $n=1,3,5$, a Taylor series expansion nearby $u_x\approx 0$ yields that $u_x \tanh^{-1} (\sqrt{1 + u_x^2})\approx 0$, which yields a well-defined energy formula. 
    However, for $n=2,4$, we obtain that $\tanh^{-1} (\sqrt{1 + u_x^2})= \infty $ for $u_x=0$, which yields an ill-defined formula. 
    These examples display the limitations of both Matano's and our methods. 
    }}
    \label{tab2}
\end{table}
\textcolor{black}{Recall that the equation \eqref{MCFpoly} is non-degenerate for bounded solutions, and therefore the characteristic equations \eqref{characMCF} do not encounter the equilibrium obstacle. On one hand, if $p_0\neq 0$ and $n$ is odd, then solutions of the characteristics are global and a singularity is not reached in finite time. On the other hand, if $p_0\leq 0$ and $n$ is even, then the characteristics display finite time blow-up, which may suppress the well-definition of a Lyapunov function, as can be seen in the Table \ref{tab2}. }
}

\subsection{\texorpdfstring{$\rho$ - } -Laplacian diffusion with an external forcing}\label{sec:rhopoly}

Consider the $\rho$-Laplacian equation \textcolor{black}{with external forcing of type $u_x^n$ with $n\in\mathbb{N}$,}
\begin{equation}\label{lap1}
    u_t = \partial_\rho u+u_x^n= (\rho -1 )|u_x|^{\rho -2}u_{xx} + u_x^n,
\end{equation}
where $\rho \geq 2 $ and $n \geq 0$.
In terms of the formulation in the previous section, we have that $ f_q=(\rho -1 )|u_x|^{\rho -2},  F^0 = - u_x^n, \text{ and } F^1 = u_t $.
%

Hence the characteristic equations \eqref{quasishootlyap} are given by
\begin{equation}\label{pLapPoly}
    \begin{aligned}
    \dot{x}&=(\rho -1 )\, |p|^{\rho -2},\\
    \dot{u}&= (\rho -1 )\, |p|^{\rho -1},\\
    \dot{p}&=- p^n,
\end{aligned}
\end{equation}
and \eqref{Lyapg} becomes 
\begin{equation}\label{pLapPoly2}
    \dot{g}= np^{n-1}.
\end{equation}
\textcolor{black}{Note that the equation \eqref{pLapPoly} encounters the equilibrium obstacle. Indeed, whenever $p=0$, we obtain that $\dot{x}=\dot{u}=\dot{p}=\dot{g}=0$. Moreover, note that whenever $p=0$ for some $x\in [0,1]$, this amounts to $u_t=0$ for such point $x\in [0,1]$, due to \eqref{lap1}. 
On one hand, if $p_0=0$, then we consider a constant $g\equiv g_0$, due to \eqref{pLapPoly2}, which amounts to $L_{pp}=(\rho - 1)\exp(g_0)|p|^{\rho-2} $. 
On the other hand, if $p_0>0$, then one can solve the equations \eqref{pLapPoly} and find $g$ by the methods of characteristics, since $\dot{p}=-p^n$ implies that $p\to 0$ as $\tau\to\infty$, and therefore the equilibrium obstacle is not reached in finite time.
For $p_0\neq 0$, note that the relevant equations in \eqref{pLapPoly} }
%
%
%
coincide with \eqref{pLapPoly}, we obtain the same solutions in \eqref{pLapPolysol} and \eqref{gpoly}. Hence the equations \eqref{exp}, \eqref{G0} with $p_*=p_0\neq 0$, and Dirichlet boundary conditions yield 
\begin{equation}
    L_{pp} = (\rho -1) \exp(g_0) |p_0|^n \, |p|^{\rho -n-2},
    \quad L^0 = - \exp(g_0) |p_0|^n u, \quad \text{and}\quad L^1 = 0.
\end{equation}
%
Hence the Lagrangian $L$ can be obtained according to \eqref{L}, yielding the following energy formula, up to a multiplicative constant $\exp(g_0) |p_0|^n$:
\begin{equation}\label{Epoly}
    E =
    \begin{cases}
        \displaystyle{\int}_0^1 \dfrac{(\rho - 1)}{(\rho -n)(\rho - n-1)} |u_x|^{\rho - n} - u \ dx\, &\text{ for } n\neq \rho, \rho-1,\\
        \displaystyle{\int}_0^1  (\rho-1)|u_x| \left(\log|u_x|-1\right)-u \ dx \, & \text{ for } n=\rho-1,\\
        \color{black}\displaystyle{\int}_0^1 (1 - \rho) \log |u_x|-u \ dx \ &  \color{black}{\text{ for} }\  \color{black} n = \rho
    \end{cases}
\end{equation}
which decays according to
\begin{equation}\label{Edecpoly}
   \frac{dE}{dt} = -\int_0^1\dfrac{u_t^2}{|u_x|^n} \, dx.
\end{equation}
\textcolor{black}{
For $\rho=2$ and $n> 2$, see \cite{Souplet} for the construction in case of a reaction term $|u_x|^n$. Moreover, for $\rho=2$ and $n\in (0,1)$, the same energy \eqref{Epoly} with decay \eqref{Edecpoly} was obtained in \cite{Laurencot}.
See both \cite{Souplet,Laurencot} for a discussion on other values of $n$ and in case of a signed reaction $a|u_x|^{n}$ for some $a\in\mathbb{R}$.
For an interplay between a gradient and Hamiltonian reaction, i.e. $h(u,u_x)=\epsilon (u^m)_x + u^n$, see \cite{FilaSacks96}.
For $\rho=2$, the equation \eqref{lap1} is not degenerate and thus the Lyapunov function regular, as solutions of a strict parabolic equation are also regular for $t>0$. However, for $\rho\neq 2$, instead of directly obtaining a Lyapunov function which is intrinsic for the degenerate equations, the authors in \cite{Attouchi,Stinner} resourced to viscosity approximations and an associated approximating Lyapunov function.} \textcolor{black}{See \cite{CIL} for a user's guide on viscosity solutions.} 

Since we have not proved any further regularity 
of the energy $E$, its derivative is also formal.
For equilibria, $u_t\equiv 0$, the energy $E$ in \eqref{Epoly} is constant, 
and its derivative $dE/dt$ given by \eqref{Edecpoly} either vanishes or it attains the value $-\infty$ in case the integrand in \eqref{Edecpoly} is not integrable, which means the derivative is not well defined. 
Similarly for time dependent solutions: either \eqref{Edecpoly} is integrable yielding negative bounded values, or \eqref{Edecpoly} is not integrable and thereby ill-defined.
Thus, equilibria may be critical points of a non-differentiable energy. 
\textcolor{black}{See \cite{Laurencot}, who mentions that \eqref{Edecpoly} is singular and it is not clear how to give a meaning to it. However, \cite[Proposition 9]{Laurencot} provides a weaker result which is sufficient to obtain dynamical information.}
Indeed, note that one can still obtain dynamical information for continuous Lyapunov functions, 
see \cite[Chapter 4]{Henry81} and \cite[Chapter 5.7]{BabinVishik92}. See also \cite{Souplet,Laurencot,Stinner}. 

Note that the energy in this example remains true for $n<0$, even though the reaction term in the vector field is singular when $u_x=0$. In this case, the decay rate of the energy in \eqref{Edecpoly} is bounded along bounded solutions of \eqref{lap1}. Thus, our methods can be applied in certain cases of quenching phenomena, whenever the hypothesis \eqref{par} holds true and one can solve the characteristic equations. See \cite{VazquezGalaktionov01} for an example of quenching in a fully nonlinear equation.

\textcolor{black}{Note that our construction can be replicated for more general external forcing. For example, when $\rho=2$ and the nonlinearity is of exponential type, see \cite{ZhangLi}.}

\subsection{Porous Medium equation}\label{sec:PME}
Consider the porous medium equation (PME) for $m\geq 1$,
\begin{equation}\label{PME}
    u_t =(u^m)_{xx} =  mu^{m-1}u_{xx}+m(m-1)u^{m-2}u_x^2.
\end{equation}
Note this is a degenerate parabolic equation for non-negative solutions $u\geq 0$, only.\footnote{The diffusion given by $(|u|^{m-1}u)_{xx}=m|u|^{m-1}u_{xx}+m(m-1)|u|^{m-3}u u_x^2$ is a natural extension that takes sign-changing solutions into account, which is thereby called the \emph{signed PME} in \cite{Vazquez06}. For the sake of simplicity, we proceed with the \textcolor{black}{non-signed} PME in the main text.}
\textcolor{black}{Instead of considering $(u^m)_{xx}$ as a nonlinear diffusion operator, we split it into two separate terms: a nonlinear degenerate diffusion, $mu^{m-1}u_{xx}$, and a nonlinear reaction, $m(m-1)u^{m-2}u_x^2$. Indeed,}
in the previous setting, $f_q = mu^{m-1} $, $F^0 = -m(m-1)u^{m-2}u_x^2 $ and $F^1 = u_t$. Therefore the characteristic equations \eqref{quasishootlyap} are given by
\begin{equation}\label{PMEchar}
    \begin{aligned}
           \dot{x} &= mu^{m-1}, \\
           \dot{u} &= mu^{m-1}p, \\
           \dot{p} &= -m(m-1)u^{m-2}p^2,
    \end{aligned}
\end{equation}
and \eqref{Lyapg} reduces to
\begin{equation}\label{PMEg}
    \dot{g} = m(m-1)u^{m-2}p.
\end{equation}
\textcolor{black}{Note that the equation \eqref{PMEchar} encounters the equilibrium obstacle. Indeed, whenever $u=0$, we obtain that $\dot{x}=\dot{u}=\dot{p}=\dot{g}=0$. Moreover, note that whenever $u=0$ for some $x\in [0,1]$, this amounts to $u_t=0$ for such point $x\in [0,1]$, due to \eqref{PME}. 
On one hand, if $u_0=0$ for $\tau\in (\tau_m,\tau_M)$, then we consider a constant $g\equiv g_0$ for $\tau\in (\tau_m,\tau_M)$, due to \eqref{PMEg}, which amounts to $L_{pp}=m\exp(g_0)u^{m-1}$. 
On the other hand, if $u_0>0$ for $\tau\in (\tau_m,\tau_M)$, then one can solve the equations \eqref{pLapPoly} and find $g$ by the methods of characteristics.}
%
For $u_0\neq 0$, we introduce the variable $\tilde{\tau}$ such that $d\tilde{\tau}/d\tau=mu^{m-2}$, with notation $(.)'=d(.)/d\tilde{\tau}$, the characteristic equations become 
\begin{equation}
    \begin{aligned}
           x' &= u, \\
           u' &= up, \\
           p'&=-(m-1)p^2, 
    \end{aligned}
\end{equation}
and
\begin{equation}\label{g'PME}
    g'= (m-1)p.
\end{equation}
We can solve this explicitly, which yields
\begin{equation}\label{pPME}
    p(\tilde{\tau}) =\frac{1}{1/p_0 + (m-1)\tilde{\tau}},
\end{equation}
and
\begin{equation}\label{gPME}
    g(\tilde{\tau}) =g_0+\log\left((m-1)p_0 \tilde{\tau}+1\right).
\end{equation}
Hence
\begin{equation}
    g(p)=g_0+\log\left(\frac{p_0}{p}\right).
\end{equation}
Thus equations \eqref{exp}, \eqref{G0} with $p_*=0$, and the Dirichlet boundary imply that 
\begin{equation}
    L_{pp} = m \exp(g_0) \left|\frac{p_0}{p}\right| u^{m-1}, \quad L^0 = 0,
    \quad \text{and}\quad L^1 = 0.
\end{equation}
%
%
Hence the Lagrangian $L$ can be obtained according to \eqref{L}, yielding the following energy candidate, up to a multiplicative constant $|p_0| \exp(g_0)$,
\begin{equation}\label{PMELyap}
    E =  \int_0^1 m u^{m-1} |u_x|\left(\log |u_x| - 1\right) \, dx,
\end{equation}
which decays according to
\begin{equation}\label{EdecPME}
   \frac{dE}{dt} = - \int_0^1\dfrac{u_t^2}{|u_x|} \, dx.
\end{equation}
Note that \eqref{PMELyap} is different from the usual energy given by $\tilde{E}=\int_0^1 \frac{u^{m+1}}{m+1} \, dx$ such that $\frac{d\tilde{E}}{dt} = -\int_0^1 [(u^{m})_x]^2\, dx$.
%
The energy $E$ decays with respect to the $L^2$-norm of $u_t$ with weight $1/|u_x|$, whereas $\tilde{E}$ decays with respect to the $L^2$-norm of $(u^m)_x$. 
Recall that the decay rate \eqref{EdecPME} is a formal computation and the energy $E$ may not be differentiable. However, one is still able to infer dynamical properties from a continuous Lyapunov function; see the discussion after equation \eqref{Edecpoly}. 
Thus, the new energy $E$ in \eqref{PMELyap} may be more suitable than $\tilde{E}$ to infer dynamical properties of the porous medium equation, such as \cite{AronsonCrandall82,DiazDiaz,Vazquez06}\textcolor{black}{, especially in case of further gradient-dependent forcing. }

{\color{black}
Yet another Lyapunov function of a rescaled porous medium equation was found in \cite{Vazquez04} for $m>1$, along with the asymptotic classification of solutions in one dimension. Indeed, considering the following rescaling, $u(t,x)=t^{-\frac{1}{m-1}}\theta(\tau,x)$ where $ t=\exp(\tau)$, the function $\theta$ satisfies the following PDE:
\begin{equation}
    \theta_\tau =(\theta^m)_{xx}+\frac{1}{m-1}\theta.
\end{equation}        
This equation possess the following Lyapunov function,
\begin{equation}
    V:=\int_0^1  \frac{|(|\theta^{m-1}|\theta)_x|^2}{2} -\frac{m}{(m+1)(m-1)} |\theta|^{m+1} dx,
\end{equation}
which decays according to 
\begin{equation}
    \frac{dV}{dt}=-m\int_0^1 |\theta|^{m-1} (\theta_\tau)^2 dx \leq 0 .
\end{equation}
The Lyapunov function $V$ decays in a similar manner as $E$, i.e., the energy $V$ decreases except for $\theta \equiv 0$ and equilibria $\theta_\tau=0$.

For an equivalence between the porous medium equation and the $\rho$-laplacian, see \cite{Iagar}. For the doubly nonlinear equation, which combines the diffusion of the porous medium and the $\rho$-laplacian, see \cite{EstebanVazquez86,Kalashnikov87}. In particular, since we have obtained a new energy for the porous medium equation, we also expect to obtain a new formula for the doubly nonlinear equation and generalizations thereof.}

\vspace{0.6cm}

\textbf{Acknowledgments.} 
We are grateful for the suggestions of Jia-Yuan Dai that significantly improved our paper.
PL was firstly supported by FAPESP, 17/07882-0, and later on by Marie Skłodowska-Curie Actions, UNA4CAREER H2020 Cofund, 847635, 
with the project DYNCOSMOS. EB was firstly supported by CNPq, 135896/2020-7, and later on by FAPESP, 23/07941-7.

\textbf{Data Availability Statement.}  We confirm that all relevant data are included in the article.

\medskip


\end{document}